# Counting Descents in Standard Young Tableaux

Ilia Barahovski

Department of Mathematics, Bar-Ilan University, Israel

email: barahilia@gmail.com

Abstract

This paper deals with the distribution of descent number in standard Young tableaux of certain shapes. A simple explicit formula is presented for the number of tableaux of any shape with two rows, with any specified number of descents. For shapes of three rows with one cell in the third row recursive formulas are given, and are solved in certain cases.

## 1. Introduction

1.1. Descent functions on standard Young tableaux were the subject of recent works [1], [6]. Expected value and variance were computed for these functions; in particular, for two important functions: the descent number and the major index. In the case of major index, the well known Stanley hook formula [9] gives an explicit generating function. However, the distribution of the descent number is unknown [1]. In this paper we find the distribution of descent number for certain shapes: hook, two rows, and three rows with one cell in the third row. Some of the results are explicit, and some are in form of recursive formulas.

1.2. The paper present extract from an MSc thesis, prepared under the supervision of Prof. Ron Adin. The current version contains mostly results without complete proofs.

## 1. Preliminaries

2.1. Let $n$ be a positive number. A **partition** of $n$ is a vector $\lambda = (\lambda_1, \lambda_2, ..., \lambda_k)$ of positive integers, where $\lambda_1 \geq \lambda_2 \geq ... \geq \lambda_k$ and $\lambda_1 + ... + \lambda_k = n$. A **Young diagram** of shape $\lambda$ is a subset of a rectangular table such that the $i$ th row contains $\lambda_i$ cells, with all rows starting at the first column. A **Young tableau** of shape $\lambda$ is obtained by inserting the integers $1, 2, ..., n$ as entries in the cells of the Young diagram of shape $\lambda$ without repetitions. The Young tableau is called **standard** if its entries increase along all rows and columns. This paper deals with standard tableaux only, so if not written, standard implied.

In the following example, $T_1$ is a Young tableau that is not standard, while $T_2$ is standard:

$T_1$
| 2 | 3 | 5 | 9 |
|---|---|---|---|
| 8 | 1 | 4 |   |
| 7 |   |   |   |
| 6 |   |   |   |

$T_2$
| 1 | 3 | 6 | 7 |
|---|---|---|---|
| 2 | 4 | 8 |   |
| 5 |   |   |   |
| 9 |   |   |   |

The **hook length** of cell $(i,j)$ is the number: $h_{i,j} = \lambda_i + \lambda'_j - i - j + 1$. In this formula $\lambda_i$ is the length of the $i$'th row and $\lambda'_j$ denotes the length of the $j$'th column.

The famous Frame-Robinson-Thrall formula [8, theorem 3.1.2] gives the number of all standard Young tableaux of shape $\lambda$ (denoted by $f^\lambda$) via the hook lengths:

$$f^\lambda = \frac{n!}{\prod_{(i,j)\in\lambda} h_{i,j}}$$

2.2. Let $T$ be the standard Young tableau. Entry $i$ is called **descent**, if $(i+1)$ placed in one of rows below $i$. $D(T)$ denotes the set of all descents in $T$. The **descent number** $des(T)$ is the number of members in $D(T)$ and **major index** $maj(T)$ is their sum:

$$des(T) = \#D(T), \qquad maj(T) = \sum_{i \in D(T)} i.$$

The following Stanley hook formula [9, corollary 7.21.5] describes the distribution of major index of standard Young tableaux:

$$\sum_{shape(T)=\lambda} q^{maj(T)} = q^{\sum_i (i-1)\lambda_i} \cdot \frac{\prod_{k=1}^{n} [k]_q}{\prod_{(i,j)\in\lambda} [h_{i,j}]_q},$$

where for each $m$: $[m]_q := 1 + q + q^2 + \ldots + q^{m-1} = \frac{q^m - 1}{q - 1}$.

There are formulas for the expected value and the variance of the descent number over all standard Young tableaux of a given shape [1, theorem 1.1], [6, theorem 2.1]. However, its exact distribution is unknown. The purpose of this research was to approach this problem.

We denote by $R_\lambda^d(condition)$ the number of standard Young tableaux of shape $\lambda$ with $d$ descents that satisfy an additional condition.

3. Results

3.1. For some simple shapes the problem may be solved using simple combinatorial observations. To start consider tableaux of hook shape as in the picture:

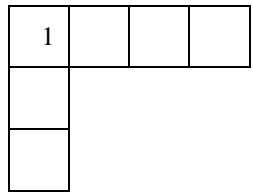

More formally, the shape of tableau is $\lambda = (n, 1^m) = (n, \underbrace{1, \ldots, 1}_{m})$, where $m \geq 0$, $n \geq 1$. The $m$ entries in the first column, below the first row, must be an increasing subsequence of $2, \ldots, n+m$ in order for the tableau to be standard. These entries completely determine the tableau; thus there are $\binom{n+m-1}{m}$ standard tableaux of shape $\lambda$. For each entry $i > 1$ in the column, $(i-1)$ is a descent. Thus each tableau of shape $\lambda$ has descent number $m$. Special cases are shapes $\lambda = (n)$ of one row or $\lambda = (1^m)$ of one column.

3.2. The next case is the class of two rows. Here the solution is more difficult. First of all, we consider the subclass of tableaux with two rows of the same length.

Denote by $R^d_{(n,n)}(k)$ the number of standard Young tableaux of shape $(n,n)$ with $d$ descents, such that the number at the end of the first row is $k$. Clearly, $R^d_{(n,n)}(k) = 0$ when the parameters $(n,k,d)$ are not in the range:

$$\begin{cases} n \geq 1, \\ 1 \leq d \leq n, \\ n \leq k < 2n. \end{cases}$$

We define an operation of **expansion** of a tableau of shape $(n-1, n-1)$ by the addition of one column to the end of the tableau. In this way we obtain a recursive formula for the distribution of descents in tableaux of shape $(n,n)$:

$$R^d_{(n,n)}(k) = R^d_{(n-1,n-1)}(k-1) + \sum_{1 \leq a < k-1} R^{d-1}_{(n-1,n-1)}(a),$$

holding for $n \geq 2$, with initial value: $R^1_{(1,1)}(1) = 1$.

The recursive formula leads to a closed form formula for this distribution:

$$R^d_{(n,n)}(k) = \begin{cases} \dfrac{2n-k}{d-1}\binom{k-n-1}{d-2}\binom{n-1}{d-2}, & \text{if } d \geq 2, k \geq n+d-1; \\ 1, & \text{if } d = 1, k = n; \\ 0, & \text{otherwise.} \end{cases}$$

This formula was discovered using an affine transformation of the parameters $(n,k,d)$, where for small values we made use of the online encyclopedia for integer sequences [6]. Summing over $k$ we finally obtain ($1 \leq d \leq n$):

$$R^d_{(n,n)} = \sum_k R^d_{(n,n)}(k) = \frac{1}{d}\binom{n-1}{d-1}\binom{n}{d-1}.$$

This formula gives the distribution of the descent number in all tableaux of two equal-length rows. These are the Narayana numbers [10].

3.3. We advanced to the class of Young tableaux of shape $\lambda = (n,m)$ with two rows of arbitrary lengths, $n \geq m \geq 1$. The key idea for the solution was to find a combinatorial connection to tableaux with two rows of the same length, which were already investigated. It turns out that adding large numbers to the end of the second row (that is, the shortest one) indeed gives tableaux with the desired shape. Moreover, the change in descent number depends only on the number at the end of the first row, and this operation is invertible. The connection may be stated as:

$$R^d_{(n,m)}(k) = \begin{cases} R^d_{(n,n)}(k), & \text{if } k < n+m; \\ R^{d+1}_{(n,n)}(k), & \text{if } k = n+m, \end{cases}$$

where $k$ once again denotes the entry at the end of the first row.

This makes it possible to find the distribution of the descent number for all shapes $(n,m)$ (when $1 \leq d \leq m \leq n$):

$$R_{(n,m)}^d = \frac{n-m+1}{d}\binom{m-1}{d-1}\binom{n}{d-1},$$

and the aim of this work is achieved for all shapes of two rows.

3.4. The treatment of tableaux with three rows is much more complicated. The class studied consists of tableaux of three rows, with only one cell in the third row. We use $R_{(n,m,1)}^d(k,c)$ to denote the number of standard Young tableaux of shape $(n,m,1)$ with $d$ descents such that $k$ is in the last cell of the first row and $c$ is in the single cell in the third row. When the number $c$ in the last row is maximal, it is possible to remove it and get a tableau with two rows (studied before) and one descent less, so: $R_{(n,m,1)}^d(k,c) = R_{(n,m)}^{d-1}(k)$. When $c$ is larger than $k$, the number at the end of the first row, it is possible to remove numbers from the second row until $c$ will be maximal and reduce to the previous case. Once again the number of descents changes by a simple law. If $k = c-1$ we get:

$$R_{(n,m,1)}^d(k,c) = R_{(n,n)}^d(k) = \frac{2n-k}{d-1}\binom{k-n-1}{d-2}\binom{n-1}{d-2},$$

and if $k < c-1$:

$$R_{(n,m,1)}^d(k,c) = R_{(n,n)}^{d-1}(k) = \frac{2n-k}{d-2}\binom{k-n-1}{d-3}\binom{n-1}{d-3}.$$

These combinatorial operations cannot be used when $k > c$.

3.5. Another approach to the class of tableaux with three rows with a single cell in the last row is again recursion. For the classes $(n,n,1)$ and $(n,m,1)$, when $n \geq m \geq 1$, it is possible to add one or more cells to the end of the tableaux and get standard Young tableaux in the same class. The change in the number of descents depends on the numbers at the ends of rows. In this way, two more recursive formulas were constructed, one for each of the two classes $(n,n,1)$ and $(n,m,1)$.

It's easy to see that the following simple conditions are necessary for $R_{(n,m,1)}^d(k,c) \neq 0$:

$$\begin{cases} n \geq m \geq 1 \\ n \leq k \leq n+m+1 \\ 3 \leq c \leq n+m+1 \\ 2 \leq d \leq m+1 \\ c \neq k \end{cases}$$

If in addition $m < n$, then $R_{(n,m,1)}^d(k,c)$ could be calculated from the following recursion:

$$R_{(n,m,1)}^d(k,c) = \begin{cases} R_{(n-1,m,1)}^d(k-1,c-1) + \sum_{a<k-1} R_{(n-1,m,1)}^{d-1}(a,c-1), & \text{if } k < c-1; \\ \sum_{a<k} R_{(n-1,m,1)}^d(a,c-1), & \text{if } k = c-1; \\ R_{(n-1,m,1)}^d(k-1,c) + \sum_{a<k-1} R_{(n-1,m,1)}^{d-1}(a,c), & \text{if } c < k < n+m+1; \\ \sum_a R_{(n-1,m,1)}^d(a,c), & \text{if } k = n+m+1. \end{cases}$$

The initial conditions are for $m = n$, where a more complicated recursion holds for $n \geq 2$:

$$R^d_{(n,n,1)}(k,c) = \begin{cases} \sum_{a<k} R^{d-1}_{(n-1,n-1,1)}(a,c), & k = 2n; \\ \sum_{a<k} R^d_{(n-1,n-1,1)}(a,c-1), & c = k+1; \\ \sum_{a<k} R^{d-1}_{(n-1,n-1,1)}(a,c-2), & k = 2n-1, c = 2n+1; \\ R^d_{(n-1,n-1,1)}(k-1,c-2) + \sum_{a<k-1} R^{d-1}_{(n-1,n-1,1)}(a,c-2), & c = 2n+1, k \notin \{2n-1, 2n\}; \\ R^d_{(n-1,n-1,1)}(k-1,c-1) + \sum_{a<k-1} R^{d-1}_{(n-1,n-1,1)}(a,c-1), & k+1 < c < 2n+1; \\ R^d_{(n-1,n-1,1)}(k-1,c) + \sum_{a<k-1} R^{d-1}_{(n-1,n-1,1)}(a,c), & c < k < 2n. \end{cases}$$

The initial value (for $n = 1$) is: $R^2_{(1,1,1)}(1,3) = 1$.

We have managed to solve these recursive formulas only in the cases that were already studied using a combinatorial approach (see section 2.4 above). The general case is still open.